\documentclass[reqno,10pt]{amsart}
\usepackage{latexsym,amssymb,amsfonts,amsmath,amsthm}
\usepackage[usenames]{color}
\usepackage{graphicx}

\newtheorem{thm}{Theorem}[section]

\newtheorem{lem}[thm]{Lemma}

\theoremstyle{definition}

\theoremstyle{remark}
\newtheorem{rem}[thm]{Remark}
\numberwithin{equation}{section}
\newcommand{\Real}{\mathbb R}

\newcommand{\F}{\mathcal{F}}

\renewcommand{\P}{\mathsf{P}}

\newcommand{\E}{\mathsf{E}}

\newcommand{\Borel}{\mathcal{B}}

\begin{document}

\title[]{Distribution of the Brownian motion on its way to hitting zero}%
\author{P. Chigansky}%
\address{Department of Statistics,
The Hebrew University,
Mount Scopus, Jerusalem 91905,
Israel}
\email{pchiga@mscc.huji.ac.il}

\author{F. Klebaner}
\address{ School of Mathematical Sciences,
Monash University Vic 3800,
Australia}
\email{fima.klebaner@sci.monash.edu.au}

\thanks{Research supported by the Australian Research Council Grant  *DP0881011}


\date{17 October,  2008}%
%\dedicatory{}%
%\commby{}%
% ----------------------------------------------------------------
\begin{abstract}
For the one-dimensional Brownian motion $B=(B_t)_{t\ge 0}$, started at $x>0$, and the  first hitting time
$\tau=\inf\{t\ge 0:B_t=0\}$, we find the probability density of $B_{u\tau}$ for a $u\in(0,1)$, i.e. of
 the Brownian motion on its way to hitting zero.
\end{abstract}
\maketitle

\section{Introduction}
The following problem has been recently addressed in \cite{JKS07a},
\cite{JKS07b}. The authors considered a continuous time subcritical
branching process $Z=(Z_t)_{t\ge 0}$, starting from the  initial
population of size $Z_0=x$. As is well known, $Z_t$ gets extinct at
the random time $T=\inf\{t\ge 0:Z_t=0\}$, and $T<\infty$ with
probability one. What can be said about $Z_{T/2}$, i.e. the
population size on the half-way to its extinction? While the
complete characterization of the law of $Z_{uT}$ with $u=1/2$, or
more generally $u\in(0,1)$, does not seem to be trackable, it turns
out that under quite general conditions
\begin{equation}
\label{JKS}
x^{u-1} Z_{uT} \xrightarrow[x\to \infty]{d} C c^{-u}e^{-u\eta},
\end{equation}
where the convergence is in distribution, $C$ and $c$ are constants,
explicitly computable in terms of parameters of $Z$ and $\eta$ is a
random variable with Gumbel distribution.

In this note we study the analogous problem for one-dimensional
Brownian motion $B=(B_t)_{t\ge 0}$, started from $x>0$. Hereafter we
assume that $B$ is defined on the canonical probability space
$(\Omega,\F,\P_x)$ and  let $\tau$ denote the first time it hits
zero, i.e. $\tau=\inf\{t\ge 0: B_t=0\}$.

\begin{thm}\label{thm}
For $x>0$ and $u\in(0,1)$, the distribution $\P_x\big(B_{u\tau}\le y\big)$ is absolutely continuous with
the density
\begin{equation}
\label{p}
p(u,x;y)=\frac{
4\sqrt{u(1-u)}xy^2
}{\pi\big\{(y-x)^2(1-u)+y^2 u\big\}\big\{(y+x)^2(1-u)+y^2 u\big\}}.
\end{equation}
\end{thm}

\begin{rem}
Notice that $p(u,x;y)$ decays as $\propto 1/y^2$ and hence its mean is infinite.
Such behavior, of course, stems from the possibility of large excursions of $B$ from the origin, before hitting zero.

\end{rem}

%\begin{figure}
%  % Requires \usepackage{graphicx}
%  \includegraphics[scale=0.5]{puxy.eps}\\
%  \caption{$p(u,x;y)$ for $x=1$ and several values of $u$\label{fig-1}}
%\end{figure}

\begin{rem}
The formula \eqref{p} implies that $x^{-1} B_{u\tau}$ has the same law under $\P_x$ as $B_{u\tau}$ under $\P_1$,
or using different notations,
\begin{equation}\label{scaling}
x^{-1}B^x_{u\tau(x)}\stackrel{d} {=} B^1_{u\tau(1)},
\end{equation}
where $B^x$ stands for the Brownian motion, starting at $x>0$, and $\tau(x)=\inf\{t\ge 0:B^x_t=0\}$
(i.e. $B^1_{u\tau(1)}$ has the density $p(u,1,y)$).
This scale invariance does not seem to be obvious at the outset and
should be compared to \eqref{JKS}, where the  scaling depends on $u$
and holds only in the limit.
\end{rem}

In the following section  we shall give an elementary proof of
Theorem \ref{thm}. In Section \ref{sec-3} our result is discussed in the
context of Doob's $h$-transform conditioning.

\section{Proof}
Let $\delta>0$ and define\footnote{$\lfloor x\rfloor$ stands for the integer part of $x\in\Real$ and
$\lceil x\rceil:=\lfloor x \rfloor +1$}
$\tau_\delta := \delta\lfloor \tau/\delta\rfloor$.
Recall that $\tau$ has the  probability density (see e.g. \cite{BoSa96}):
\begin{equation}
\label{tauf}
f(x;t)=\frac {\partial}{\partial t}\P_x(T\le t)=\frac{x}{\sqrt{2\pi t^3}}e^{-x^2/2t},\quad t\ge 0, \quad x>0.
\end{equation}
Let $\hat M_{s,t}:=\inf_{s\le r<  t}B_r$ and $\phi(\cdot)$ be a  continuous bounded
function, then\footnote{$I(\cdot)$ denotes the indicator function}
\begin{align*}
\E_x \phi\big(B_{u\tau_\delta}\big)&=\sum_{k= 0}^\infty \E_x \phi\big(B_{u\tau_\delta}\big)
I\Big(\tau\in \big[\delta k, \delta(k+1)\big)\Big)\\
&=\sum_{k= 0}^\infty \E_x \phi\big(B_{u\delta k}\big)
I\Big(\tau\in \big[\delta k, \delta(k+1)\big)\Big)\\
&=
\sum_{k= 0}^\infty \E_x \phi\big(B_{u\delta k}\big)
I\Big(\hat M_{0,\delta k}>0, \hat M_{\delta k, \delta (k+1)}\le 0\Big)
\\
&=\sum_{k= 0}^\infty \E_x \phi\big(B_{u\delta k}\big)
I\Big(\hat M_{0,u\delta k}>0\Big)\P_x\Big( \hat M_{u\delta k,\delta k}>0, \hat M_{\delta k, \delta (k+1)}\le 0\big|
\F^B_{u\delta k}
\Big)\\
&=\sum_{k= 0}^\infty \E_x \phi\big(B_{u\delta k}\big)
I\Big(\hat M_{0,u\delta k}>0\Big)\P_x\Big( \hat M_{u\delta k,\delta k}>0, \hat M_{\delta k, \delta (k+1)}\le 0\big|
B_{u\delta k}
\Big)\\
&=\phi\big(x\big)\P_{x}\Big(\tau\in [0, \delta )\Big)
+
\sum_{k= 1}^\infty \E_x \phi\big(B_{u\delta k}\big)
I\Big(\hat M_{0,u\delta k}>0\Big)\times \\
&\hskip 2.0in
\P_{B_{u\delta k}}\Big(\tau\in \big[(1-u)\delta k, (1-u)\delta k +\delta\big)\Big)
\end{align*}
\begin{align}\label{longeq}
&=\phi\big(x\big)\int_{0}^{\delta } f\big(x;t\big)dt+
\sum_{k= 1}^\infty \E_x \phi\big(B_{u\delta k}\big)
I\Big(\hat M_{0,u\delta k}>0\Big)\int_{(1-u)\delta k}^{(1-u)\delta k + \delta} f\big(B_{u\delta k};t\big)dt
\nonumber
\\
&=\phi\big(x\big)\int_{0}^{\delta } f\big(x;t\big)dt+
\int_0^\infty\phi(y) \bigg\{\sum_{k= 1}^\infty\int_{\delta k}^{\delta (k+1)} q(x,u\delta k,y) f\big(y;t-u\delta k\big)dt\bigg\}dy
\nonumber
\\
&=\phi\big(x\big)\int_{0}^{\delta } f\big(x;t\big)dt+
\int_0^\infty\phi(y)\bigg\{ \int_{\delta}^{\infty} q(x,u \lfloor t/\delta\rfloor \delta,y) f\big(y;t-u\lfloor t/\delta\rfloor \delta\big)dt\bigg\}dy
\nonumber
\\
&=\phi\big(x\big)\int_{0}^{\delta } f\big(x;t\big)dt+
\int_0^\infty\phi(y)\bigg\{ \int_{0}^{\infty} q(x,u \lceil t/\delta\rceil\delta,y) f\big(y;t+\delta-u\lceil t/\delta\rceil\delta\big)dt\bigg\}dy,
\end{align}
where $q(x,t,y)$ is the probability density of $\P_x\big(\hat M_{0,t}>0,B_t\in dy \big)$ with respect to the
Lebesgue measure (see e.g. formula 1.2.8 page 126, \cite{BoSa96}):
\begin{equation}
\label{qxty}
q(x,t,y)=
\bigg\{
\frac{1}{\sqrt{2\pi t}}e^{-(y-x)^2/2t}-\frac{1}{\sqrt{2\pi t}}e^{-(y+x)^2/2t}\bigg\},
\quad x,y>0.
\end{equation}
By continuity of the densities \eqref{tauf} and \eqref{qxty}, for any fixed $x>0$ and $u\in(0,1)$, the function
$$
F_\delta(t,y):=q(x,u \lceil t/\delta\rceil\delta,y) f\big(y;t+\delta-u\lceil t/\delta\rceil\delta\big)
$$
converges to
$$
\lim_{\delta \to 0}F_\delta(t,y)=q(x,u t,y) f(y;t-ut), \quad \forall t\ge  0,\ \ y\ge 0.
$$
In  Lemma \ref{lem} below we exhibit a function $G(t,y)$, independent of $\delta$,
such that
\begin{equation}\label{domdom}
\F_\delta(t,y)\le G(t,y), \quad \forall (t,y)\in \Real^2_+\quad \text{and}\quad
\int_{\Real^2_+} G(t,y)dtdy<\infty,
\end{equation}
and hence, the dominated convergence and \eqref{longeq} imply
\begin{equation}\label{dc}
\lim_{\delta \to 0}\E_x \phi\big(B_{u\tau_\delta}\big) =
\lim_{\delta\to 0}\int_{\Real^2_+}\phi(y)F_\delta(t,y)dydt =
\int_{\Real^2_+}\phi(y) q(x,u t,y) f(y;t-ut)dtdy.
\end{equation}
On the other hand,  $\lim_{\delta\to 0}\tau_\delta= \tau$, $\P_x$-a.s. and thus  by continuity of $B_t$,
$\lim_{\delta\to 0}B_{u\tau_\delta}=B_{u\tau}$, $\P_x$-a.s. for any $u\in(0,1)$. Thus, by arbitrariness of $\phi$,
\eqref{dc}  implies that the distribution of $B_{u\tau}$ has the density:
$$
p(u,x;y):=\int_0^\infty q(x,u t,y) f(y;t-ut)dt.
$$
A calculation now yields:
\begin{align*}
p(u,x;y)&=
\int_0^\infty
\frac{y}{\sqrt{2\pi} \big(t(1-u)\big)^{3/2}}
e^{-y^2/2t(1-u)}
\bigg\{\frac{1}{\sqrt{2\pi ut}}e^{-(y-x)^2/2ut}-\frac{1}{\sqrt{2\pi ut}}e^{-(y+x)^2/2ut}\bigg\}dt
\\
&=\frac{y}{2\pi (1-u)^{3/2}u^{1/2}}
\int_0^\infty
\frac{1}{t^2}
\bigg\{e^{-(y-x)^2/2ut-y^2/2t(1-u)}-e^{-(y+x)^2/2ut-y^2/2t(1-u)}\bigg\}dt,
\end{align*}
and by a change of variables
\begin{align*}
p(u,x;y)&=
\frac{y}{2\pi (1-u)^{3/2}u^{1/2}}
\bigg\{
\Big(\frac{(y-x)^2}{2u}+\frac{y^2}{2(1-u)}\Big)^{-1}-
\Big(\frac{(y+x)^2}{2u}+\frac{y^2}{2(1-u)}\Big)^{-1}\bigg\}\\
&=\frac{2y\sqrt{u}}{2\pi \sqrt{1-u}}
\bigg\{\frac{1}{(y-x)^2(1-u)+y^2 u}-
\frac{1}{(y+x)^2(1-u)+y^2 u}
\bigg\}\\
&=\frac{
\sqrt{u(1-u)}4xy^2
}{\pi\Big\{(y-x)^2(1-u)+y^2 u\Big\}\Big\{(y+x)^2(1-u)+y^2 u\Big\}}.
\end{align*}

The statement of the Theorem \ref{thm} now follows from:

\begin{lem}\label{lem}
\eqref{domdom} holds with $G(t,x)$ defined in \eqref{bnd} below.
\end{lem}
\begin{proof}
Set $t^\delta :=\lceil t/\delta\rceil\delta$, so that
$t\le t^\delta\le t+\delta$, and
\begin{align*}
F_\delta(t,y)&\le
\frac{1}{\sqrt{ u t^\delta}}\bigg\{
e^{-(y-x)^2/2u t^\delta}-e^{-(y+x)^2/2u t^\delta}\bigg\}
\frac{y}{ (t+\delta-ut^\delta)^{3/2}}e^{-\frac 1 2 y^2/(t+\delta-ut^\delta)}
\\
&\le I(t\le \delta)\frac{1}{\sqrt{ u \delta}}e^{-(y-x)^2/(2u \delta)}
\frac{y}{ \big((1-u)\delta\big)^{3/2}}e^{-\frac 1 2 y^2/(\delta+(1-u)\delta)}+
\\
&\hskip 0.15in I(t> \delta)
\frac{1}{\sqrt{ u t}}
\bigg\{
e^{-(y-x)^2/2u t^\delta}-e^{-(y+x)^2/2u t^\delta}\bigg\}\frac{y}{ \big((t+\delta)(1-u)\big)^{3/2}}e^{-\frac 1 2 y^2/(t(1-u)+\delta)}\\
&=: I(t\le \delta)A+I(t>\delta)B.
\end{align*}
Since the function $z^2e^{-Cz}$ with $C>0$ attains its maximum
$4e^{-2}/C^2$ on the interval $[0,\infty)$ at $z:=2/C$,
\begin{multline*}
A=
\frac{y}{\sqrt{ u (1-u)^3}}
\frac{1}{ \delta^2}
\exp\bigg\{-\bigg(\frac{(y-x)^2}{2u }+\frac 1 2 \frac{y^2}{(2-u)}\bigg)\frac 1\delta\bigg\}\le \\
\frac{y}{\sqrt{ u (1-u)^3}}
\bigg(\frac{(y-x)^2}{2u }+\frac 1 2 \frac{y^2}{(2-u)}\bigg)^{-2}.
\end{multline*}
Similarly, for $t> \delta$,
\begin{align*}
B&=
\frac{y}{\sqrt{u(1-u)^3t(t+\delta)^3}}
e^{-(y-x)^2/2u t^\delta}\bigg\{
1-e^{-2xy/u t^\delta}\bigg\}
e^{-\frac 1 2 y^2/(t(1-u)+\delta)}  \\
&\le\frac{y}{\sqrt{u(1-u)^3}}\frac{1}{ t^2}
e^{-(y-x)^2/2u (t+\delta)}\bigg\{
1-e^{-2xy/u t}\bigg\}
e^{-\frac 1 2 y^2/(t(1-u)+\delta)} \\
& \le\frac{y}{\sqrt{u(1-u)^3}}\frac{1}{ t^2}
e^{-(y-x)^2/4u t}\bigg\{
1-e^{-2xy/u t}\bigg\}
e^{-\frac 1 2 y^2/(t(2-u))}.
\end{align*}
Hence  for $\delta\in(0,1]$ we have the bound
\begin{multline}\label{bnd}
F_\delta(t,y)\le
\frac{y}{\sqrt{ u (1-u)^3}}
\bigg(\frac{(y-x)^2}{2u }+\frac 1 2 \frac{y^2}{(2-u)}\bigg)^{-2}I(t\le 1)
+\\
\frac{y}{\sqrt{u(1-u)^3}}\frac{1}{ t^2}
e^{-(y-x)^2/4u t}\bigg\{
1-e^{-2xy/u t}\bigg\}
e^{-\frac 1 2 y^2/(t(2-u))}=:G(t,y).
\end{multline}
Since for $u\in(0,1)$ and $x>0$, the quadratic function is lower bounded:
$$
\frac{(y-x)^2}{u}+ \frac{y^2}{(2-u)}\ge  \frac{x^2}{2},
$$
the first function in the right hand side of \eqref{bnd} is
integrable on $\Real_+^2$.
Further,
\begin{align*}
\int_0^\infty &\frac{y}{\sqrt{u(1-u)^3}}\frac{1}{ t^2}
e^{-(y-x)^2/4u t}\bigg\{
1-e^{-2xy/u t}\bigg\}
e^{-\frac 1 2 y^2/(t(2-u))}dt\\
&=\frac{y}{\sqrt{u(1-u)^3}}
\bigg\{
\bigg(
\frac{(y-x)^2}{4u}+
\frac{ y^2}{2(2-u)}
\bigg)^{-1}-
\bigg(
\frac{(y-x)^2}{4u}+
\frac{ y^2}{2(2-u)}+
\frac{2xy}{u}
\bigg)^{-1}
\bigg\}
\\
&=
\frac{4yu(2-u)}{\sqrt{u(1-u)^3}}
\bigg\{
\frac{1}{(y-x)^2(2-u)+2y^2u}
-
\frac{1}{(y-x)^2(2-u)+2y^2u+8xy(2-u)}
\bigg\}
\\
&=\frac{32u(2-u)^2 xy^2/\sqrt{u(1-u)^3}}{\Big\{(y-x)^2(2-u)+2y^2u\Big\}\Big\{(y-x)^2(2-u)+2y^2u+8xy(2-u)\Big\}}.
\end{align*}
The latter function decays as $\propto 1/y^2$ as $y\to\infty$ and is bounded away from zero, uniformly in $y\ge 0$,
and thus is integrable on $\Real_+$. Since the last term in the right hand side of \eqref{bnd} is nonnegative,
by Fubini theorem it is an integrable function on $\Real_+^2$ for all $u\in (0,1)$ and $x>0$.
\end{proof}

\section{A connection to Doob's $h$-transform}\label{sec-3}

In this section we show that the random variable $B_{u\tau}$ has the same
density as the so called {\em scaled Brownian excursion} at the corresponding time, averaged over its length.
The latter process is defined by conditioning in the sense of Doob's $h$-transform, and it would be natural
to identify this formal conditioning with the usual conditional probability.
While in the analogous discrete time setting, such identification is evident, its precise justification
in our case remains an open problem.

For a fixed time $T>0$, let $R=(R_t)_{t\le T}$ be the 3-dimensional Bessel bridge $R=(R_t)_{t\le T}$, starting at $R_0=x$ and
ending at zero. Namely, $R$ is the radial part\footnote{$\|\cdot\|$ denotes the Euclidian norm in $\Real^n$}
\begin{equation}\label{R}
R_t = \|V_t\|, \quad t\in[0,T],
\end{equation}
of the 3-dimensional Brownian bridge $V=(V_t)_{t\le T}$ with $V_0=v$ and $V_T=0$:
$$
V_t = v+W_t-\frac t T(W_T+v), \quad t\in [0,T],
$$
where  $v\in\mathbb{R}^3$ with $\|v\|=x$ and  $W$ is a standard Brownian motion in
$\mathbb{R}^3$.

The law of $R$ coincides with the law of the scaled Brownian excursion process, which is defined
as ``the Brownian motion, started at $x>0$ and conditioned to hit zero for the
first time at time $T$''. Here the conditioning is understood in the sense of Doob's $h$-transform
(see Ch. IV, \S 39, \cite{RW}, and \cite{Bl92}, \cite{Do84} for the in depth treatment).

On the other hand, one can speak on the regular conditional measure induced on the space of Brownian excursions
(started from $x>0$), given $\tau=\inf\{t\ge 0: B_t=0\}$. More precisely,
let $E$ be a subset of continuous functions $C_{[0,\infty)}(\Real)$, such that for all $\omega\in E$, $\omega(0)=x$
and for each $\omega$ there is a positive number $\ell(\omega)$, called the {\em excursion length},
such that $\omega(t)>0$ for $0<t< \ell(\omega)$ and  $\omega(t)\equiv 0$ for all $t\ge \ell(\omega)$. $E$ together with the smallest $\sigma$-algebra $\mathcal{E}$, making all
coordinate mappings measurable, is called the {\em excursion} space  (see \S 3, \cite{BoSa96} for the
brief reference and \cite{Bl92} for  more details). Let $\mu_x(T,\cdot)$, $T\in[0,\infty)$ be a probability kernel on the
excursion space $(E,\mathcal{E})$, i.e. a family of measures such that $T\mapsto \mu_x(T,A)$ is a measurable function for all $A\in \mathcal{E}$
and $\mu_x(T,\cdot)$ is a probability measure on $\mathcal{E}$ for each $T\ge 0$. By definition, $\mu_x(T,\cdot)$ is
a regular conditional probability of $B_{t\wedge \tau}$ given $\tau$, if for any bounded and measurable functional $F$
on $(E,\mathcal{E})$:
$$
\E_x F(B_{\cdot \wedge \tau})I(\tau\in A) = \int_A \int_E F(\omega) \mu_x(s,d\omega) f(x;s)ds,\quad \forall  A\in \Borel(\Real),
$$
where $f(x;t)$ is the density of $\tau$, defined in \eqref{tauf}. In particular, for any bounded measurable
function $\phi$ and some $u\in(0,1)$,
\begin{equation}
\label{mumu}
\E_x \phi(B_{u\tau}) = \int_0^\infty \int_E \phi\big(\omega(us)\big) \mu_x(s,d\omega) f(x;s)ds.
\end{equation}
We were not able to trace any general result, from which the identification  of $\mu_x(T,d\omega)$ with the probability
$\nu_x(T,d\omega)$, induced on $(E,\mathcal{E})$ by the aforementioned Bessel bridge $R$, could be deduced.
While the latter, of course, is intuitively appealing, its precise justification remains elusive
(some relevant results can be found in \cite{FPY93}). The calculations below show that
\begin{equation}
\label{nunu}
\int_0^\infty \int_E \phi\big(\omega(us)\big) \nu_x(s,d\omega) f(x;s)ds=\E_x \phi(B_{u\tau}),
\end{equation}
indicating in favor of such identification.

For a fixed $T>0$ and $u\in(0,1)$, the distribution of $R_{uT}$, i.e. the restriction of $\nu_x(T,d\omega)$ to the time $t:=uT$,
has a density $q_{uT}(x;y)$ with respect to the Lebesgue measure $dy$, which can be computed as follows.
We have
$$
E V_t = v(1-t/T), \quad \mathrm{cov}(V_t)=I\frac{t(T-t)}{T},
$$
where $I$ is $3$-by-$3$ identity matrix.
Notice that the law of the Bessel bridge $R$ in \eqref{R} doesn't depend on the particular $v$ as
long as $\|v\|=x$ and
it will be particularly convenient to carry out the calculations
for the specific choice  $v=(x,0,0)$. Fix a constant $u\in(0,1)$  and let $\xi_1,\xi_2,\xi_3$ be i.i.d.
standard Gaussian random variables.
Then, for $t:=uT$,
$$
R_{uT} \stackrel{d}{=} \sqrt{\Big(\sqrt{Tu(1-u)}\xi_1+x(1-u)\Big)^2+
Tu(1-u)\xi_2^2+Tu(1-u)\xi_3^2
}.
$$
The random variable $\theta:=\xi_2^2+\xi_3^2$ has $\chi_2^2$
distribution, which is the same as the exponential distribution with
parameter $1/2$ and hence
\begin{equation}
\label{Rx}
R_{uT} \stackrel{d}{=} b\sqrt{\Big(\xi+a\Big)^2+
\theta
},
\end{equation}
where $\xi$ is written for $\xi_1$ and $a:=x\sqrt{(1-u)/Tu}$ and $b:=\sqrt{(1-u)Tu}$ are defined for brevity.
The density of $(\xi+a)^2$ is given by:
\begin{multline*}
f_1(z):=\frac{d}{dz}P\Big((\xi+a)^2\le z\Big)=\frac{d}{dz}
\int_{-\sqrt{z}}^{\sqrt{z}}\frac 1 {\sqrt{2\pi}} e^{-(x-a)^2/2}dx= \\
\frac 1 {2\sqrt{2\pi}\sqrt{z}} \Big(e^{-(\sqrt{z}-a)^2/2}
+ e^{-(\sqrt{z}+a)^2/2}\Big).
\end{multline*}
The density  of $(\xi+a)^2+\theta$ is the convolution of $f_1$ and the exponential density with parameter $1/2$:
\begin{align*}
f_3(y)&:=\int_0^y f_1(z)f_2(y-z)dz=
\int_0^y \frac 1 {2\sqrt{2\pi}\sqrt{z}} \Big(e^{-(\sqrt{z}-a)^2/2}
+ e^{-(\sqrt{z}+a)^2/2}\Big)\frac 1 2 e^{-1/2 (y-z)}dz\\
&=\int_0^{\sqrt{y}} \frac 1 {2\sqrt{2\pi}} \Big(e^{-(\tilde z-a)^2/2}
+ e^{-(\tilde z+a)^2/2}\Big) e^{-1/2 (y-\tilde z^2)}d\tilde z\\
&=\frac {e^{-a^2/2- y/2}} {2\sqrt{2\pi}}\int_0^{\sqrt{y}}  \Big(e^{\tilde z a }
+ e^{-\tilde z a}\Big) d\tilde z=
\frac {e^{-a^2/2- y/2}} {\sqrt{2\pi}}\int_0^{\sqrt{y}}  \cosh(\tilde za) d\tilde z\\
&=\frac {e^{-a^2/2- y/2}} {\sqrt{2\pi} a}\sinh(\sqrt{y}a).
\end{align*}
Consequently, the density of $\sqrt{(\xi+a)^2+\theta}$ is given by
$$
f_4(z) := 2z f_3(z^2)=
\frac{\sqrt{2}e^{-a^2/2}}{\sqrt{\pi} a}ze^{- z^2/2}\sinh(za),
$$
and, finally, the density  of $b\sqrt{(\xi+a)^2+\theta}$
is
$$
f_5(z):=\frac 1 b f_4(z/b)=\frac{\sqrt{2}e^{-a^2/2}}{\sqrt{\pi} a b^2}ze^{- z^2/2b^2}\sinh(za/b)=
\frac{z}{\sqrt{2\pi} ab^2}\bigg\{e^{-(a-z/b)^2/2}-e^{-(a+z/b)^2/2}\bigg\}.
$$
Hence by \eqref{Rx}, $R_{uT}$ has the density
\begin{multline*}
q_{uT}(x;y)=\frac{y}{\sqrt{2\pi} x(1-u)}\frac{1}{\sqrt{Tu(1-u)}}
\bigg\{\exp\bigg(-\frac{1}{T}\frac{\big(x(1-u)-y\big)^2}{2u(1-u)}\bigg)
-\\
\exp\bigg(-\frac{1}{T}\frac{\big(x(1-u)+y\big)^2}{2u(1-u)}\bigg)\bigg\},
\end{multline*}
We shall abbreviate by writing
$$
C_1=\frac{\big(x(1-u)-y\big)^2}{2u(1-u)}, \ C_2=\frac{\big(x(1-u)+y\big)^2}{2u(1-u)},\
C_3=\frac{y}{\sqrt{2\pi} xu^{1/2}(1-u)^{3/2}}.
$$
Then
\begin{multline*}
\int_0^\infty \int_E \phi\big(\omega(ut)\big) \nu_x(t,d\omega) f(x;t)dt=
\int_0^\infty \phi(y)\int_0^\infty q_{ut}(x;y)f(x;t)dt\\
=\int_0^\infty \phi(y)\int_0^\infty
C_3\frac{1}{t^{1/2}}
\Big\{e^{-C_1/t}-e^{-C_2/t}\Big\}
\frac{x}{\sqrt{2\pi} t^{3/2}}
e^{-x^2/2t}dtdy
=:\int_0^\infty \phi(y)p(u,x;y)dy
\end{multline*}
and by a change of  variables,
\begin{multline*}
p(u,x;y)=C_3\frac{x}{\sqrt{2\pi}}
\int_0^\infty
\Big\{e^{-(C_1+x^2/2)t}-e^{-(C_2+x^2/2)t}\Big\}
t^{-2}dt= \\
C_3\frac{x}{\sqrt{2\pi}}
\bigg\{\frac{1}{C_1+x^2/2}-\frac{1}{C_2+x^2/2}\bigg\}.
\end{multline*}
A calculation yields,
$$
C_1+x^2/2=
\frac{
(1-u)(x-y)^2+y^2u}{2u(1-u)}\quad \text{and}\quad
C_2+x^2/2=\frac{
(1-u)(x+y)^2+y^2u}{2u(1-u)},
$$
and, consequently,
\begin{align*}
p(u,x;y)&=
\frac{y}{\sqrt{2\pi} xu^{1/2}(1-u)^{3/2}}
\frac{x}{\sqrt{2\pi}}
\bigg\{
\frac{2u(1-u)}{
(1-u)(x-y)^2+y^2u}
-
\frac{2u(1-u)}{
(1-u)(x+y)^2+y^2u}
\bigg\}
\\
&=\frac{4xy^2\sqrt{u(1-u)}}{\pi\big\{(1-u)(x-y)^2+y^2u\big\}\big\{(1-u)(x+y)^2+y^2u\big\}},
\end{align*}
which in view of \eqref{mumu} and  \eqref{p}, imply  \eqref{nunu}.
% ----------------------------------------------------------------

%\bibliographystyle{plain}
%\bibliography{fimabm}

\end{document}